\documentclass[10pt]{amsart}
\usepackage{latexsym}
\usepackage{amssymb}
\usepackage{amsmath}
\usepackage{mathrsfs}
\usepackage{xcolor}
\usepackage{graphicx}
\usepackage{bbm}
\usepackage{bbold}
\def\epi{\mathop{\fam0 epi}\nolimits}
\def\cop{\mathop{\fam0 cop}\nolimits}
\def\dom{\mathop{\fam0 dom}\nolimits}
\def\On{\mathop{\fam0 On}\nolimits}

\font\new=labi1200

\begin{document}

\title{Convexity and Cone-Vexing
}

\author{S.~S. Kutateladze}
\begin{abstract}
This is a talk delivered on September 20, 2007   at the
conference ``Mathematics in the Modern World'' on the occasion of the
fiftieth anniversary of the Sobolev Institute of Mathematics in Novosibirsk, Russia.
\end{abstract}
\address[]{
Sobolev Institute of Mathematics\newline
\indent Novosibirsk, 630090
\indent RUSSIA}
\email{
sskut@member.ams.org
}

\maketitle

\begin{itemize}
\item[]{{\large\bf\color{blue} To Vex}\  {\scshape\color{red} (WordWeb 5.0)}}

\item[\bf1.]{
 Cause annoyance in; disturb, especially by minor irritations}

\item[\bf2.]{
 Disturb the peace of mind of; afflict with mental agitation or distress}

\item[\bf3.]{
Change the arrangement or position of}

\item[\bf4.]{
 Subject to prolonged examination, discussion, or deliberation\\
   ``vex the subject of the death penalty"}

\item[\bf5.]{
 Be a mystery or bewildering to\\
   ``a vexing problem''}

\end{itemize}

\medskip

\section{\large\bf\color{blue} Agenda}

Convexity stems  from the remote ages
and reigns in geometry, optimization, and functional analysis.
The union of abstraction and convexity has produced
abstract convexity which is a vast area of today's research,
sometimes profitable but sometimes bizarre. Cone-vexing is
a popular fixation of vexing  conic icons.

The idea of convexity feeds generation, separation, calculus,
and approximation.  Generation appears as duality; separation, as optimality;
calculus,  as representation; and approximation, as stability.
This is an overview of the origin, evolution,
and trends  of convexity.

Study of convexity in the Sobolev Institute was initiated by
Leonid Kantorovich (1912--1986) and Alexandr Alexandrov (1912--1999).
This talk is a part of their memory.

\section{\large\bf\color{blue} {\new Elements}, Book I}
Mathematics resembles linguistics sometimes and pays tribute to etymology,
hence, history.  Today's convexity is a centenarian, and abstract convexity
is much younger.

Vivid convexity is  full of abstraction, but
traces back to the idea of a solid figure which stems from Euclid.
Book~I of his {\it Elements}~\cite{Euclid} has expounded plane geometry and
defined a~boundary and a~figure as follows:

\begin{itemize}
\item[]{\small
Definition 13. A boundary is that which is an extremity of anything.}

\item[]{\small
Definition 14. A figure is that which is contained by any boundary or boundaries.}
\end{itemize}

\section{\large\bf\color{blue}{\new Elements}, Book XI}
Narrating solid geometry in Book~XI, Euclid travelled
in the opposite direction from  solid to  surface:

\begin{itemize}
\item[]{\small
Definition 1. A solid is that which has length, breadth, and depth.}

\item[]{\small
Definition 2. An extremity of a solid is a surface}.
\end{itemize}

He  proceeded with the relations of similarity and equality for solids:

\begin{itemize}
\item[]{\small
Definition 9. Similar solid figures  {\it are  those contained
by similar planes equal in multitude}}.

\item[]{\small
Definition 10. Equal and similar solid figures {\it are those contained by
similar planes equal in multitude and magnitude.}}
\end{itemize}

\section{\large\bf\color{blue} The Origin of Convexity}
Euclid's definitions  seem vague, obscure, and even unreasonable
if applied to the figures other than convex polygons and polyhedra.
Euclid also introduced a formal concept of ``cone'' which has
a well-known natural origin. However, convexity was ubiquitous
in his geometry by default. The term ``conic sections''  was coined
as long ago as 200 BCE by  Apollonius of Perga. However, it was long
before him that Plato had formulated his famous allegory of cave~\cite{Plato}.
The shadows on the wall  are often convex.

Euclid's  definitions imply the intersection of half-spaces. However, the concept
of intersection belongs to set theory which appeared only at the  end
of the nineteenth century.
It is wiser to seek for the origins
of the ideas of Euclid in his past rather than his future.
Euclid was a scientist not a foreteller.

\section{\large\bf\color{blue} Harpedonaptae}
The predecessors of Euclid are the  harpedonaptae of Egypt as often
sounds   at the lectures on the history of mathematics.
The harpedonaptae or rope-stretchers measured tracts of land
in the capacity of surveyors. They administered {\it cadastral surveying\/}
which gave rise to the notion of geometry.
If anyone stretches a rope that surrounds however many  stakes, he will
distinguish a convex polygon, which is  up to infinitesimals
a~typical compact convex set or abstract subdifferential of the present-day
mathematics. The rope-stretchers discovered convexity experimentally  by measurement.
Hence,   a few words are in order about these forefathers
of their Hahn--Banach next of kin of today.

\section{\large\bf\color{blue} The History of Herodotus}

Herodotus   wrote in Item 109 of Book II {\it Enerpre}~\cite{Herodotus} as follows:
\begin{itemize}
\item[]{\small
Egypt was cut up: and they said that this king
distributed the land to all the Egyptians, giving an equal square
portion to each man, and from this he made his revenue, having
appointed them to pay a certain rent every year: and if the river
should take away anything from any man's portion, he would come to the
king and declare that which had happened, and the king used to send
men to examine and to find out by measurement how much less the piece
of land had become, in order that for the future the man might pay
less, in proportion to the rent appointed: and I think that thus the
art of geometry was found out and afterwards came into Hellas also.}
\end{itemize}

\section{\large\bf\color{blue} Sulva Sutras}
Datta~\cite{Datta}  wrote:

\begin{itemize}
\item[]{\small
\dots One who was well versed in that science was called in
ancient India as samkhyajna (the expert of numbers), parimanajna
(the expert in measuring), sama-sutra-niranchaka
(uniform-rope-stretcher), \linebreak Shul\-ba-vid (the expert in Shulba) and
Shulba-pari\-prcchaka (the inquirer into the Shulba).}
\end{itemize}

Shulba also written as \'Sulva or Sulva was in fact the geometry of vedic times
as codified in $\text{\it \'S}ulva\ S\bar{u}tras$.

\section{\large\bf\color{blue} Veda}
Since ``veda'' means  knowledge, the vedic epoch and literature are indispensable
for understanding the origin and rise of mathematics.
In 1978  Seidenberg~\cite{Seidenberg} wrote:

\begin{itemize}
\item[]{\small
Old-Babylonia [1700 BC] got the theorem of Py\-thagoras from India or that both
Old-Babylonia  and India got it from a third source. Now the Sanskrit scholars
do not give me a date  so far back as 1700 B.C. Therefore I postulate
a pre-Old-Babylonian (i.e., pre-1700 B.C.) source of the kind of geometric
rituals we see preserved in the Sulvasutras, or at least for
the mathematics involved in these rituals.}
\end{itemize}

Some recent facts and evidence prompt us that the roots of rope-stretching
spread in a~much deeper past than we were accustomed to acknowledge.

\section{\large\bf\color{blue} Vedic Epoch}
The exact chronology still evades us and
  Kak~\cite{Kak} commented  on the Seidenberg paper:

\begin{itemize}
\item[]{\small
That was before archaeological finds disproved
the earlier assumption of a break in Indian civilization in the second millennium B.C.E.;
it was this assumption of the Sanskritists that led Seidenberg to postulate a third earlier
source. Now with our new knowledge, Seidenberg's conclusion of India being the source
of the geometric and mathematical knowledge of the ancient world fits in with the new
chronology of the texts.}

\item[]{\small
\dots in the absence of conclusive evidence, it is prudent to take the
most conservative  of these dates, namely 2000 B.C.E. as the latest
period to be associated with the Rigveda.}
\end{itemize}

\section{\large\bf\color{blue} Mathesis and Abstraction}
Once upon  a time mathematics was everything. It is not now but
still carries the genome of {\it mathesis universalis}.
Abstraction is the mo\-ther of reason
and the gist of mathematics.
It enables us to collect the particular
instances  of any many with some property we observe or study.
Abstraction entails generalization and proceeds by analogy
which is tricky and might be misleading. Inventory of the true
origins of
any instance of abstraction  is in order from time to time.

``Scholastic'' differs from ``scholar.'' Abstraction is limited
by taste, tradition, and common sense. The challenge of abstraction
is alike the call of freedom. But no freedom is exercised  in solitude. The holy gift of abstraction
coexists with  gratitude and respect to the legacy of
our predecessors who collected  the gems of reason
and saved them in the treasure-trove of mathematics.

\section{\large\bf\color{blue} Enter Abstract Convexity}
Stretching a rope taut between two stakes produces a closed straight line segment
which is the continuum in modern parlance.
Rope-stret\-ching raised the problem of measuring the continuum.
The continuum hypothesis of
set theory is  the shadow of the ancient problem of harpedonaptae.
Rope-stretching  independent  of  the position of stakes is
uniform with respect to direction in space.  The mental experiment
of uniform rope-stretching   yields  a~compact convex figure.
The harpedonaptae were experts in convexity.


Convexity has found solid grounds in set theory.
The Cantor paradise became an official residence of convexity.
Abstraction becomes an axiom of set theory. The abstraction axiom
enables us to reincarnate a property, in other words,  to collect and
comprehend.  The union of convexity and abstraction was inevitable.
Their child is  abstract convexity~\cite{MD}--\cite{IoRu}.

\section{\large\bf\color{blue}Minkowski Duality}

Let $\overline{E}$  be a~vector lattice
$E$ with the adjoint top $\top:=+\infty$ and bottom $\bot:=-\infty$.
Assume further that $H$  is some subset of $E$ which is by implication a~(convex)
cone in $E$, and so the bottom of $E$
lies beyond~$H$. A subset $U$  of~$H$ is {\it convex relative to~}
$H$ or $H$-{\it convex\/}
provided that $U$ is the $H$-{\it support set\/}
$U^H_p:=\{h\in H:h\le p\}$ of some element $p$ of $\overline{E}$.

Alongside the $H$-convex sets we consider
the so-called $H$-convex elements. An element   $p\in \overline{E}$
is  $H$-{\it convex} provided that $p=\sup U^H_p$;~i.e., $p$
represents the supremum of the $H$-support set of~$p$.
The $H$-convex elements comprise the cone which is denoted by
$\mathscr C(H,\overline{E}$).  We may omit  the references to $H$ when $H$ is clear
from the context. It is worth noting that
convex elements and sets are ``glued together''
by the {\it Minkowski duality\/} $ \varphi:p\mapsto U^H_p$.
This duality enables us to study convex elements and sets simultaneously
\cite{Fenchel}.

\section{\large\bf\color{blue}Enter the Reals}
Optimization is the science of choosing the best.
To choose, we use preferences. To optimize, we use infima and
suprema (for bounded subsets) which is practically
the {\it least upper bound property}.
So optimization needs ordered sets and primarily
Dedekind  complete lattices.

To operate with preferences,
we use group structure. To aggregate and scale, we use linear structure.

All these are happily provided by the {\it reals\/} $\mathbb R$, a one-dimensional
Dedekind complete vector lattice. A~Dedekind  complete vector
lattice is a~{\it Kantorovich space}.

\section{\large\bf\color{blue}Legendre in Disguise}
An abstract minimization problem  is as follows:
$$
x\in X,\quad f(x)\rightarrow \inf.
\eqno{(*)}
$$
Here $X$ is a vector space and
$f: X\rightarrow \overline{\mathbb R}$ is a numeric
function taking possibly infinite  values.
The sociological trick  includes the problem
into a parametric family  yielding
the  {\it Young--Fenchel transform\/} of $f$:
$$
f^*(l):=\sup_{x\in X}{(l(x)-f(x))},
$$
of $l\in X^{\#}$,  a linear functional over $X$.
The epigraph of $f^*$  is a convex subset of
$ X^{\#}$ and so $f^*$  is  convex.
Observe that $-f^*(0)$ is the value of $(*)$.

\section{\large\bf\color{blue}Order Omnipresent}
A convex function is locally a positively homogeneous convex
function, a {\it sublinear functional}. Recall that
$p: X\to\mathbb R$ is sublinear whenever
$$
\epi p:=
\{ (x,\ t)\in X\times\mathbb R: p(x)\le t\}
$$
is a cone. Recall that a numeric function is uniquely
determined from its epigraph.

Given
$C\subset X$, put
$$
H(C):=\{(x,\ t)\in X\times\mathbb R^+ : x\in tC\},
$$
the  {\it H\"ormander transform\/} of $C$~\cite{Her}.
Now, $C$ is convex if and only if $H(C)$  is a~cone.
A~space with a cone is a {\it $($pre$)$ordered vector space}.

``The order, the symmetry, the harmony enchant us\dots.''
(Leibniz)

\section{\large\bf\color{blue}Fermat's Criterion}

The {\it subdifferential\/} of $f$
at  $\bar x$ is defined as

$$
\partial f(\bar x):=\{l\in X^\# :(\forall x\in X)\
l(x)-l(\bar x)\leq f(x)-f(\bar x) \}.
$$

A point $\bar x$ is a solution
to the  minimization problem $(*)$ if and only if
$$
0\in \partial f(\bar x).
$$

This {\it Fermat criterion} turns into the Rolle Theorem in a smooth case
and is  of little avail without effective tools for calculating $\partial f(\bar x).$
A~convex analog of the  ``chain rule''  is in order.

\section{\large\bf\color{blue}Enter Hahn--Banach}
The {\it Dominated Extension}, an alias of Hahn--Banach, takes the form
$$
\partial (p\circ \iota) (0)=(\partial p)(0)\circ \iota,
$$
with $p$ a~sublinear functional over $X$ and
$\iota $ the identical embedding of some~subspace of~$X$ into~$X$.


If the target $\mathbb R$ may be replaced with an ordered vector space
$E$, then $E$ admits {\it dominated extension}.

\section{\large\bf\color{blue}Enter Kantorovich}

The matching of convexity and order
was established in two steps.

{\scshape Hahn--Banach--Kantorovich  Theorem.}
{\it Every Kantorovich space admits dominated extension of linear operators}.

This theorem proven by Kantorovich in~1935 was a first attractive result
of the theory of ordered vector spaces.

{\scshape Bonnice--Silvermann--To Theorem.}
{\it Each ordered vector space
admitting dominated extension of linear operators is a~Kantorovich space.}

\section{\large\bf\color{blue}Nonoblate Cones}
Consider cones
$K_1$
and
$K_2$
in a topological vector space
$X$
and put
$\varkappa:=(K_1,K_2)$.
Given a pair
$\varkappa$
define the correspondence
$\Phi_{\varkappa}$
from
$X^2$
into
$X$
by the formula
$$
\Phi_{\varkappa}:=\{(k_1,k_2,x)\in X^3:
x=k_1-k_2\in K_\imath \}.
$$
Clearly, $\Phi_{\varkappa}$ is a cone or, in other words,
a~conic correspondence.

The pair $\varkappa$ is {\it nonoblate\/}
whenever $\Phi_{\varkappa}$
is open at the zero. Since
$\Phi_{\varkappa}(V)=V\cap K_1-V\cap K_2$
for every
$V\subset X$,
the nonoblateness of $\varkappa$
means that
$$
\varkappa V:=(V\cap K_1-V\cap K_2)\cap(V\cap K_2-V\cap K_1)
$$
is a zero neighborhood  for every zero neighborhood~
$V\subset X$.

\section{\large\bf\color{blue}Open Correspondences}
Since $\varkappa V\subset V-V$, the nonoblateness of
$\varkappa$ is equivalent to the fact that the system of sets
$\{\varkappa V\}$ serves as a filterbase of zero neighborhoods while
$V$ ranges over some base of the same filter.

Let $\Delta_n:x\mapsto(x,\dots,x)$ be
 the embedding of
$X$
into the diagonal
$\Delta_n(X)$
of
$X^n$.
A pair of cones
$\varkappa:=(K_1,K_2)$
is nonoblate if and only if
$\lambda:=(K_1\times K_2,\Delta_2(X))$
is nonoblate in~$X^2$.

Cones
$K_1$
and
$K_2$
constitute a nonoblate pair if and only if the conic correspondence~
$\Phi\subset X\times X^2$
defined as
$$
\Phi:=\{(h,x_1,x_2)\in X\times X^2 :
x_\imath+h\in K_\imath\ (\imath:=1,2)\}
$$
is open at the zero.


\section{\large\bf\color{blue}General Position of Cones}
Cones~
$K_1$
and~
$K_2$
in a topological vector space~
$X$
are {\it in general position\/}
iff

{\bf (1)}~
the algebraic span of $K_1$
and~
$K_2$
is some subspace
$X_0\subset X$;
i.e.,
$X_0=K_1-K_2=K_2-K_1$;

{\bf (2)}~the subspace
$X_0$
is complemented; i.e., there exists a continuous projection
$P:X\rightarrow X$
such that
$P(X)=X_0$;

{\bf (3)}~$K_1$
and~
$K_2$
constitute a nonoblate pair in~
$X_0$.

\section{\large\bf\color{blue}General Position of Operators}
Let $\sigma_n$
stand for the  rearrangement of coordinates
$$
\sigma_n:((x_1,y_1),\dots, (x_n,y_n))\mapsto ((x_1,\dots,x_n),\linebreak
(y_1,\dots,y_n))
$$
which establishes an isomorphism between
$(X\times Y)^n$
and
$X^n\times Y^n$.

Sublinear operators  $P_1,\dots,P_n:X\rightarrow E\cup \{+\infty\}$
are {\it in general position\/} if so are
the cones $\Delta_n(X)\times E^n$
and
$\sigma_n(\epi (P_1)\times\dots\times\epi (P_n))$.

Given a cone $K\subset X$, put
$$
\pi_E(K):=\{T\in\mathscr L(X,E): Tk\leq 0\ (k\in K)\}.
$$
Clearly,
$\pi_E(K)$
is a cone in
$\mathscr L(X,E)$.

{\scshape Theorem.} {\it Let
$K_1,\dots,K_n$
be cones in a~topological vector space~$X$
and let
$E$
be a topological Kantorovich space.  If
$K_1,\dots,K_n$
are in general position then}
$$
\pi_E(K_1\cap\dots\cap K_n)=\pi_E(K_1)+\dots+\pi_E(K_n).
$$
This formula opens  a way to various separation results.

\section{\large\bf\color{blue}Separation}

{\scshape Sandwich Theorem.} {\it Let
$P,Q:X\rightarrow E\cup \{+\infty\} $
be sublinear operators in general position.
If
$P(x)+Q(x)\geq 0$
for all
$x\in X$
then there exists a~continuous linear operator~
$T:X\rightarrow E$
such that}
$$
-Q(x)\leq Tx\leq P(x)\quad  (x\in X).
$$

Many efforts were made to abstract these results to
a more general algebraic setting and, primarily,
to semigroups and semimodules. Tropicality
chases separation~\cite{IFA, CGQ}.

\section{\large\bf\color{blue}Canonical Operator}

Consider a~Kantorovich space
$E$
and an arbitrary nonempty set
$\mathfrak A$.
Denote by
$l_\infty (\mathfrak A,E)$
the set of all order bounded mappings from
$\mathfrak A$
into $E$; i.e.,
$f\in l_\infty (\mathfrak A,E)$
if and only if
$f:\mathfrak A \to E$
and
$\{f(\alpha):\alpha\in\mathfrak A\}$
is order bounded in $E$.
It is easy to verify that
$l_\infty (\mathfrak A,E)$ becomes a Kantorovich
space if endowed with the coordinatewise algebraic operations and order.
The operator
$\varepsilon_{\mathfrak A, E}$
acting from
$l_\infty (\mathfrak A,E)$
into
$E$
by the rule
$$
\varepsilon_{\mathfrak A, E}:f\mapsto\sup \{f(\alpha):\alpha\in\mathfrak A\}
\quad (f\in l_\infty (\mathfrak A,E))
$$
is called the {\it canonical sublinear operator\/}
given
$\mathfrak A$
and
$E$.
We often write
$\varepsilon_{\mathfrak A}$
instead of
$\varepsilon_{\mathfrak A, E}$
when it is clear from the context what Kantorovich space is meant.
The notation
$\varepsilon_n$
is used when the cardinality of~$\mathfrak A$
equals
$n$ and we call the operator
$\varepsilon_n$
{\it finitely-generated}.

\section{\large\bf\color{blue}Support Hull}

Consider a~set~
$\mathfrak A$
of linear operators acting from
a~vector space
$X$
into a~Kantorovich space
$E$.
The set
$\mathfrak A$
is {\it weakly order  bounded\/} if
$\{\alpha x:\alpha\in\mathfrak A\}$
is order bounded for every
$x\in X$.
Denote by
$\langle\mathfrak A\rangle x$
the mapping that assigns the element
$\alpha x\in E$
to each
$\alpha\in \mathfrak A$,
i.e.
$\langle\mathfrak A\rangle x: \alpha\mapsto\alpha x$.
If
$\mathfrak A$
is weakly order bounded then
$\langle\mathfrak A\rangle x\in l_\infty (\mathfrak A,E)$
for every fixed
$x\in X$.
Consequently, we obtain the linear operator
$\langle\mathfrak A\rangle:X\rightarrow l_\infty (\mathfrak A,E)$
that acts as
$\langle\mathfrak A\rangle:x\mapsto\langle\mathfrak A\rangle x$.
Associate with
$\mathfrak A$
one more operator
$$
p_{\mathfrak A}: x\mapsto\sup \{\alpha x: \alpha\in\mathfrak A\}\quad(x\in X).
$$
The operator
$p_{\mathfrak A}$
is sublinear. The support set
$\partial p_{\mathfrak A}$
is denoted by
$\cop (\mathfrak A)$
and referred to as  the {\it support hull\/} of
$\mathfrak A$.

\section{\large\bf\color{blue}  Hahn--Banach  in Disguise}

{\scshape Theorem.} {\it If
$p$
is a~sublinear operator with
$\partial p=\cop (\mathfrak A)$
then $
P=\varepsilon_{\mathfrak A}\circ \langle\mathfrak A\rangle.
$
Assume further that
$p_1: X\to E$
is a~sublinear operator and
$p_2: E\to F$
is an increasing sublinear operator. Then
$$
\partial (p_2\circ p_1)
=
\{ T\circ\langle\partial p_1\rangle:
 T\in L^+(l_{\infty}(\partial p_1, E),F)
\ \&\ T\circ\Delta_{\partial p_1}\in \partial p_2
\}.
$$
Moreover, if
$\partial p_1=\cop (\mathfrak A_1)$
and
$\partial p_2=\cop (\mathfrak A_2)$
then}
$$
\partial (p_2\circ p_1)
=\bigl\{T\circ\langle\mathfrak A_1\rangle : T\in L^+
(l_{\infty}(\mathfrak A_1,E),F)\
\&\
\left(\exists\alpha\in\partial\varepsilon_{\mathfrak A_2}\bigr)\
T\circ\Delta_{\mathfrak A_1}=\alpha\circ\langle\mathfrak A_2\rangle\right\}.
$$


\section{\large\bf\color{blue}Enter Boole}
Cohen's final solution of the problem of the cardinality of the
continuum within ZFC gave rise to the Boolean-valued models
by Vop\v enka, Scott, and Solovay.
Takeuti coined the term ``Boolean-valued analysis''
for applications of the new models
to functional analysis~\cite{BA}.

Let
$B$
be a~complete Boolean algebra. Given an ordinal
$\alpha$,
put
$$
V_{\alpha}^{(B)}
:=\{x:
(\exists \beta\in\alpha)\ x:\dom (x)\rightarrow
B\ \&\ \dom (x)\subset V_{\beta}^{(B)}  \}.
$$
The {\it Boolean-valued universe\/}
${\mathbb V}^{(B)}$
is
$$
{\mathbb V}^{(B)}:=\bigcup\limits_{\alpha\in\On} V_{\alpha}^{(B)},
$$
with $\On$ the class of all ordinals. The truth
value $[\![\varphi]\!]\in B$ is assigned to each formula
$\varphi$ of ZFC relativized to ${\mathbb V}^{(B)}$.

\section{\large\bf\color{blue}~Enter Descent}

Given $\varphi$, a~formula of ZFC, and
$y$, a~subset ${\mathbb V}^{B}$; put
$A_{\varphi}:=A_{\varphi(\cdot,\ y)}:=\{x:\varphi (x,\ y)\}$.
The {\it descent\/}
$A_{\varphi}{\downarrow}$
of a~class
$A_{\varphi}$
is
$$
A_{\varphi}{\downarrow}:=\{t:t\in {\mathbb V}^{(B)} \ \&\  [\![\varphi  (t,\ y)]\!]=\mathbb 1\}.
$$
If
$t\in A_{\varphi}{\downarrow}$
then
it is said
that
{\it $t$
satisfies
$\varphi (\cdot,\ y)$
inside
${\mathbb V}^{(B)}$}.


The {\it descent\/}
$x{\downarrow}$
of an element
$x\in {\mathbb V}^{(B)}$
is defined by the rule
$$
x{\downarrow}:=\{t: t\in {\mathbb V}^{(B)}\ \&\  [\![t\in x]\!]=\mathbb 1\},
$$
i.e. $x{\downarrow}=A_{\cdot\in x}{\downarrow}$.
The class $x{\downarrow}$ is a~set.
If $x$ is a~nonempty set
inside
${\mathbb V}^{(B)}$ then
$$
(\exists z\in x{\downarrow})[\![(\exists z\in x)\ \varphi (z)]\!]
=[\![\varphi(z)]\!].
$$

\section{\large\bf\color{blue} The Reals in Disguise}

There is an object
$\mathscr R$
inside
${\mathbb V}^{(B)}$ modeling $\mathbb R$, i.e.,
$$
[\![\mathscr R\ {\text{is the reals}}\,]\!]=\mathbb 1.
$$
Let $\mathscr R{\downarrow}$ be the descend of
 the carrier $|\mathscr R|$ of the algebraic system
\hbox{$\mathscr R:=(|\mathscr R|,+,\,\cdot\,,0,1,\le)$}
inside ${\mathbb V}^{(B)}$.
Implement the descent of the structures on $|\mathscr R|$
to $\mathscr R{\downarrow}$ as follows:
$$
\gathered
x+y=z\leftrightarrow [\![x+y=z]\!]=\mathbb 1;
\\
xy=z\leftrightarrow [\![xy=z]\!]=\mathbb 1;
\\
x\le y\leftrightarrow [\![x\le y]\!]=\mathbb 1;
\\
\lambda x=y\leftrightarrow [\![\lambda^\wedge x=y]\!]=\mathbb 1
\\
(x,y,z \in\mathscr R{\downarrow},\ \lambda\in\mathbb R).
\endgathered
$$

{\scshape Gordon Theorem.}
{\it $\mathscr R{\downarrow}$
with the descended structures is a universally complete
Kantorovich space with base
$\mathscr B(\mathscr R{\downarrow})$
isomorphic to~$B$}.

\section{\large\bf\color{blue}Approximation}

Convexity of harpedonaptae was stable in the sense that no variation
of stakes within the surrounding rope can ever spoil the convexity
of the tract to be surveyed.

Study of stability in abstract convexity
is accomplished sometimes by introducing various  epsilons
in appropriate places. One of the earliest excursions in this direction
is connected with the classical  Hyers--Ulam stability theorem for
$\varepsilon$-convex functions~\cite{Almost}.
Exact calculations with epsilons and sharp estimates
are sometimes bulky and slightly mysterious.  Some alternatives
are suggested
by actual infinities, which is illustrated with the conception
of {\it infinitesimal optimality}.

\section{\large\bf\color{blue}Enter Epsilon and Monad}

Assume given a~convex operator
$f:X\to E\cup{+\infty}$
and a~point
$\overline x$
in the effective domain
$\dom(f):=\{x\in X:f(x)<+\infty\}$
of~$f$.
Given
$\varepsilon \ge 0$
in the positive cone
$E_+$
of
$E$,
by the
$\varepsilon $-{\it subdifferential\/}
of~$f$
at~$\overline x$
we mean the set
$$
\partial\, {}^\varepsilon\!f(\overline x)
:=\big\{T\in L(X,E):
(\forall x\in X)(Tx-Fx\le T\overline x -f\overline x+\varepsilon) \big\},
$$
with
$L(X,E)$
standing as usual for the space of linear operators
from~$X$
to~$E$.

Distinguish
some downward-filtered subset~$\mathscr E$ of
$E$
that is composed of positive elements.
Assuming
$E$ and~$\mathscr E$
 standard, define the {\it monad\/}
$\mu (\mathscr E)$ of $\mathscr E$ as
$\mu (\mathscr E):=\bigcap\{[0,\varepsilon ]:\varepsilon \in
{}^\circ\!\mathscr E\}$.
The members of $\mu(\mathscr E)$ are {\it positive
infinitesimals\/}  with respect to~$\mathscr E $.
As usual,
${}^\circ\!\mathscr E$
denotes the external set of~all standard members of~$E$,
the {\it standard part\/} of~$\mathscr E$.

\section{\large\bf\color{blue}Subdifferential Halo}
Assume that the monad $\mu (\mathscr E )$
is an external cone over ${}^\circ \mathbb R $ and, moreover,
$\mu (\mathscr E)\cap{}^\circ\! E=0$.
In application, $\mathscr E $ is usually
the filter of order-units of $E$.
The relation of
{\it infinite proximity\/} or
{\it infinite closeness\/}
between the members of $E$ is introduced as follows:
$$
e_1 \approx e_2 \leftrightarrow e_1 -e_2 \in\mu
(\mathscr E )\ \&\ e_2 -e_1 \in\mu (\mathscr E ).
$$

Now
$$
Df(\overline x):=\bigcap\limits_{\varepsilon \in{}^\circ \mathscr E }\,
\partial _\varepsilon f(\overline x)=
\bigcup\limits_{\varepsilon \in\mu (\mathscr E )}\,
\partial _\varepsilon f(\overline x);
$$
the {\it infinitesimal subdifferential} of
$f$ at  $\overline x$.
The elements of
$Df(\overline x)$ are
{\it infinitesimal subgradients\/}
of $f$ at~$\overline x$.

\section{\large\bf\color{blue}Exeunt Epsilon}

{\scshape Theorem.} {\it
Let $f_1:X\times Y\rightarrow E\cup +\infty$ and
$f_2:Y\times Z\rightarrow E\cup +\infty$ be convex operators.
Suppose that the convolution
$f_2\vartriangle f_1$ is infinitesimally exact at some point $(x,y,z)$; i.e.,
$
(f_2\vartriangle f_1)(x,y)\approx f_1(x,y)+f_2(y,z).
$
If, moreover, the
convex sets $\epi(f_1,Z)$ and $\epi(X,f_2)$ are
in general
position then}
$$
D(f_2\vartriangle f_1)(x,y)=
Df_2(y,z)\circ Df_1(x,y).
$$

This talk bases on the recent book \cite{Subdiff} which
covers other relevant topics.

\section{\large\bf\color{blue}Models Galore}

The essence of mathematics resides in freedom, and abstraction
is the freedom of generalization. Freedom is the loftiest ideal
and idea of man, but it  is  demanding, limited, and vexing.
So is abstraction. So are its instances in convexity.
Abstract convexity starts with repudiating the heritage of harpedonaptae,
which  is annoying and vexing but may turn out rewarding.

Freedom of set theory empowered us with the Boolean-valued models
yielding  a~lot of surprising and unforeseen visualizations
of the continuum.  Many promising
opportunities are open nowadays to modeling  the powerful habits of
reasoning and verification.

Convexity is  a topical
illustration of the wisdom and strength of mathematics, the
ever  fresh art and science of calculus.

\bibliographystyle{plain}

\begin{thebibliography}{99}

\bibitem{Euclid}
Heaves~Th.  (1956)
{\it The Thirteen Books of Euclid's Elements.} Vol.~1--3.
New York: Dover Publications.

\bibitem{Plato}
Ferrari G.~R.~F. (Ed.) (2000)
{\it Plato: The Republic}.
Berkeley: University of California (Translated by Tom Griffith).

\bibitem{Herodotus}
Macaulay~G.~C. (Tr.) (1890)
{\it The History of Herodotus}.
London and New York: Macmillan.

\bibitem{Datta}
Datta B. (1993)
{\it Ancient Hindu Geometry: The Science of the Sulba}.
New Delhi: Cosmo Publishers.

\bibitem{Seidenberg}
Seidenberg A. (1978)
The origin of mathematics.
{\it Archive for History of Exact Sciences}.
{\bf18}, 301--342.

\bibitem{Kak}
Kak S.~C. (1997)
Science in Ancient India.
In:  Sridhar~S.~R. and  Mattoo~N.~K. (eds.)
{\it Ananya: A~Portrait of India}.
New York: AIA,  399--420.

\bibitem{MD}
Kutateladze~S.~S. and Rubinov~A.~M. (1972)
Minkowski duality and its applications.
{\it Russian Math. Surveys},  {\bf27}:3, 137--191.

\bibitem{KutRub}
Kutateladze~S.~S. and Rubinov~A.~M. (1976)
{\it Minkowski Duality and Its Applications}.
Novosibirsk: Nauka Publishers [in Russian].

\bibitem{Fuch}
Fuchssteiner~B. and Lusky~W. (1981)
{\it Convex Cones.} Amsterdam: North-Holland.

\bibitem{Notions}
H\"ormander L. (1994)
{\it Notions of Convexity.} Boston:  Birkh\"auser.

\bibitem{Singer}
Singer I. (1997)
{\it Abstract Convex Analysis}.
 New York: John Wiley \& Sons.

\bibitem{PR}
Pallaschke~D. and Rolewicz~S. (1998)
{\it Foundations of Mathematical Optimization, Convex
Analysis Without Linearity}. Dordrecht: Kluwer Academic Publishers.

\bibitem{Rub}
Rubinov~A.~M. (2000)
{\it Abstract Convexity and Global Optimization}. Dordrecht:
Kluwer Academic Publishers.

\bibitem{IoRu}
Ioffe~A.~D. and Rubinov~A.~M. (2002)
Abstract convexity and nonsmooth analysis. Global aspects.
{\it  Adv. Math. Econom.}, {\bf4}, 1--23.

\bibitem{Fenchel}
Fenchel W. (1953)
{\it Convex Cones, Sets, and Functions.}
Princeton: Princeton  Univ. Press.

\bibitem{Her}
H\"{o}rmander~L. (1955)
Sur la fonction d'appui  des ensembles convexes dans une espace lokalement
convexe. {\it Ark. Mat.}, {\bf3}:2, 180--186 [in French].

\bibitem{IFA}
Litvinov G.~L., Maslov V.~P., and Shpiz~G.~B. (2001)
Idempotent functional analysis: an algebraic approach.
{\it Math. Notes}, {\bf69}:5, 758--797.

\bibitem{CGQ}
Cohen~G., Gaubert~S., and Quadrat~J.-P. (2002)
{\it Duality and Separation Theorems in Idempotent Semimodules}.
Rapporte de recherche No.~4668, Le Chesnay CEDEX: INRIA Rocquetcourt. 26~p.

\bibitem{BA}
Kusraev A.~G. and Kutateladze~S.~S. (2005)
{\it Introduction to Boolean-Valued Analysis}.
Moscow: Nauka Publishers [in Russian].

\bibitem{Almost}
Dilworth S.~J., Howard R., and   Roberts J.~W. (2006)
A~general theory of almost convex functions.
{\it Trans. Amer. Math. Soc.},
{\bf 358}:8, 3413--3445.

\bibitem{Subdiff}
Kusraev A.~G. and Kutateladze~S.~S. (2007)
{\it Subdifferential Calculus: Theory and Applications}.
Moscow: Nauka Publishers [in Russian].

\end{thebibliography}

\end{document}